\documentstyle[11pt]{article}
\begin{document}
\begin{center} 

\Large 
{\bf PRINCIPAL PAIRS FOR OSCILLATORY SECOND
ORDER LINEAR DIFFERENTIAL EQUATIONS}\footnote{

Research partially supported by grants from the Natural Sciences and
Engineering Research Council (Canada) and by Research Grant \# 199404
of the Academy of Sciences of the Czech Republic}\\

\normalsize  
\vspace{.5cm} 

\bigskip 
MARTIN E. MULDOON\\
Department of Mathematics \& Statistics\\
York University\\ 
North York, Ont. M3J 1P3, Canada\\
\bigskip 

\bigskip
FRANTI\v{S}EK NEUMAN\\
Mathematical Institute  -- Brno Branch\\
Academy of Sciences of the Czech Republic\\
\v{Z}i\v{z}kova 22\\
CZ -- 616 62 Brno, Czech Republic\\

\bigskip

\end{center} 
\begin{abstract}
Nonoscillatory second order differential equations always admit
``special'', principal solutions.  For a certain type of oscillatory
equation principal pairs of solutions were introduced by \'A. Elbert, F.
Neuman and J. Vosmansk\'y, {\em Diff. Int. Equations} {\bf 5} (1992),
945--960.  In this paper, the notion of principal pair is extended to a
wider class of oscillatory equations. Also an interesting property of
some of the principal pairs is presented that makes the notion of these
``special'' pairs more understandable. 
\end{abstract} 

\parindent=1cm 

\parskip=0.4cm 
\section{Introduction}
For nonoscillatory equations, the notion of principal solution was
introduced by W. Leighton, M. Morse and P. Hartman \cite{H}. 
Principal pairs were defined for certain types of oscillatory
equations in \cite{ENV}. 

In this paper, we extend the notion of principal pairs to a wider class
of oscillatory equations. Under certain conditions, general enough to
cover many important special equations (Bessel, Airy, etc.) we also study
the question of identifying the second member of such a pair when the
first member is given. We find that among all the solutions, the second
member is the one whose zeros approximate most closely the zeros of the
derivative of the first member of the pair.  In these cases, this
property may also serve as an alternative definition, connecting in a
certain sense the Kummer and the Pr\"ufer transformations and bringing
more light to this particular choice of ``good'' pairs of solutions of
linear second order oscillatory differential equations. 

\section{Basic facts and definitions}
Consider the equation
\begin{equation} 
y'' + q(x)y = 0 \label{q}
\end{equation}
on $[x_0,\infty)$, where $q \in C^0[x_0,\infty)$.

If equation (\ref{q}) is nonoscillatory, then a solution $y_1$ of
(\ref{q}) is called {\em principal}, in accordance with the
terminology
of Leighton, Morse and Hartman (see, e.g., \cite{H}), if 
\begin{equation}
\int^\infty y_1^{-2}(x) dx = \infty,
\end{equation}
or, equivalently, if $\lim_{x \rightarrow \infty} y_1/y_2 = 0$ for
every solution $y_2$ of (\ref{q}) linearly independent of
$y_1$.  In the nonoscillatory case such a principal solutions
always exists and is unique up to multiplication by a constant
factor.

Now, consider an {\em oscillatory} equation (\ref{q}) with two
linearly
independent solutions $y_1,y_2$ normalized by having their
Wronskian
\begin{equation}
w=w(y_1,y_2)= y_1(x)y_2'(x) - y_1'(x)y_2(x)
\end{equation}
equal to $\pm 1$, i.e.,
\begin{equation}
|w| = 1. \label{w}
\end{equation}
It was proved in \cite{ENV} that if
\begin{equation}
\lim_{x \rightarrow \infty}[y_1^2(x) + y_2^2(x)] = L \ne 0,
\label{L}
\end{equation}
for a pair of 
solutions $y_1,y_2$ of an oscillatory equation (\ref{q}),
satisfying (\ref{w}),
then  the pair $y_1,y_2$ is unique up to orthogonal
transformation:
\begin{equation}
 \left(
\begin{array}{c} y_1\\ y_2 \end{array} \right) \mapsto A \left(
\begin{array}{c} y_1\\ y_2 \end{array} \right), \;{\rm with}\;A^TA
= I,
\label{dot}
\end{equation}
i.e., if $(y_1,y_2)$ is one such pair then every such pair is of
the form 
$(ay_1 +by_2,\;cy_1 +dy_2)$ with $a^2+c^2 = 1 = b^2 +d^2,\;
ab+cd =0$.  

In \cite{ENV}, each of these pairs $(y_1,y_2)$ was called a
{\em principal} pair  of solutions of (\ref{q}). 
Some sufficient conditions for the existence of principal pairs
$(y_1,y_2)$ of solutions of (\ref{q}) with $q \in C^1[x_0,\infty)$,
were also given there, based on the fact that $v(x) =
y_1^2(x) + y_2^2(x)$ is a solution of the corresponding Appell
equation (see \cite{appell}, \cite{H}) 
\begin{equation}
v''' + 4q(x)v' + 2q'(x)v = 0 \label{a}
\end{equation}
on $[x_0,\infty)$, admitting three linearly independent solutions
$y_1^2(x), y_2^2(x)$ and $y_1(x)y_2(x)$.
Here we extend the investigation of pairs of
solutions to the cases where the $L$ in (\ref{L}) is $0$ or
$+\infty$.

We shall need the following important facts concerning the Kummer
transformation of the second order equations of the form (\ref{q}),
derived by O. Bor${\rm \stackrel{\circ}{u}}$vka, and summarized in
his monograph \cite{B}; see also \cite[Chapter 2]{FN}.
Consider an equation (\ref{q}) (either oscillatory or
nonoscillatory)  and a pair of linearly independent solutions,
$y_1$ and $y_2$.  Denote by $w = w(y_1,y_2)$ their Wronskian
\begin{equation}
w(y_1,y_2)= y_1(x)y_2'(x) - y_1'(x)y_2(x) = w \ne 0.
\end{equation}
Following O. Bor${\rm \stackrel{\circ}{u}}$vka, define the (first)
phase $\alpha$ of
(\ref{q}) with respect to $y_1,y_2$ as a continuous function
$\alpha \in C^0[0,\infty)$ satisfying the relation
\begin{equation}
\tan \alpha = y_1(x)/y_2(x),
\end{equation}
whenever this quantity is defined, i.e., for $y_2(x) \ne 0$. Then
$\alpha
\in C^3[x_0,\infty),\;\alpha'(x) \ne 0$ on $[x_0,\infty)$, and
$$ y_1(x) = \varepsilon | \alpha'(x)|^{-1/2}\sin \alpha(x), $$
 $$ y_2(x) = \varepsilon | \alpha'(x)|^{-1/2}\cos \alpha(x), $$
$$ \varepsilon = \pm\sqrt{-w \cdot\; {\rm sign}\; \alpha'},$$
$ y(x) = k_1 | \alpha'(x)|^{-1/2} \sin(\alpha(x) +  k_2)$ being a
general solution of (\ref{q}) for arbitrary constants $k_1$ and
$k_2$. 
Moreover, equation (\ref{q}) is oscillatory on $[x_0,\infty)$ (as
$x
\rightarrow \infty$) if and only if $\lim_{x \rightarrow \infty}
|\alpha (x)| = \infty$; see \cite{B}.

Each solution of (\ref{q}) is bounded if and only if $
|\alpha' (x)|^{-1}$ is bounded on $[x_0,\infty)$; see \cite{FN}.

Each solution of (\ref{q}) tends to $0$ as $x \rightarrow
\infty$
if and only if $
\lim_{x \rightarrow \infty} |\alpha' (x)|^{-1} =0$; see \cite{FN}.

In particular, we have 
$$  y_1^2(x) + y_2^2(x) = - \frac{w}{\alpha'(x)} $$
on $[x_0,\infty)$ and for $y_1$ and $y_2$ with $w(y_1,y_2) = - 1$,
we get
$$ y_1(x) = \pm | \alpha'(x)|^{-1/2}\sin \alpha(x), $$
 $$ y_2(x) = \pm | \alpha'(x)|^{-1/2}\cos \alpha(x)$$
and
$$  y_1^2(x) + y_2^2(x) = (\alpha'(x))^{-1}. $$

\section{Principal pairs of solutions --- extension}
In the previous section, the number $L$ in the limit relation
(\ref{L}) was supposed to satisfy $0 < L < \infty$. This raises the
question of extension to the cases $L =0,\;+\infty$.  In both of
these cases,  the limit relation (\ref{L}) must be supplemented by
an additional condition in order to characterize sums of squares up
to orthogonal transformation.

Consider a pair $\overline{y}_1, \overline{y}_2$ of linearly
independent solutions of (\ref{q}) with Wronskian
$w(\overline{y}_1, \overline{y}_2)$ satisfying
\begin{equation}
|w(\overline{y}_1, \overline{y}_2)| =1.
\end{equation}
We can write
$$ \overline{y}_1 = ay_1 + by_2,$$
$$ \overline{y}_2 = cy_1 + dy_2,$$
where
$$
\left| {\rm det} \left( \begin{array}{cc} a&b\\c&d \end{array}
\right) w(y_1,y_2)
\right| = |ad-bc| = 1.
$$
Hence
$$ \overline{y}_1(x) = \varepsilon |\alpha'(x)|^{-1/2}[a\sin\alpha(
x)
+b\cos\alpha(x) ],$$
$$ \overline{y}_2(x) = \varepsilon |\alpha'(x)|^{-1/2}[c\sin
\alpha(x)
+d\cos \alpha (x) ],$$
and
\begin{eqnarray}
\overline{v}(x)& =& \overline{y}_1^2(x) + \overline{y}_2^2(x)
\nonumber \\
& =&
v(x) \left[ (a^2 +c^2) \sin^2 \alpha (x) + (ab+cd)
\sin 2\alpha(x) + (b^2 +d^2) \cos^2 \alpha(x)
\right]. \label{star}
\end{eqnarray}
In case $v(x) \rightarrow L$ with $0 < L < +\infty$, it is clear
that $\overline{v}(x) $ approaches a finite limit if
and only if 
$a^2 +c^2 = b^2 +d^2 = K > 0$, and $ab+cd=0$. 
However, in case $L = 0$ the condition 
$v(x) \rightarrow 0$ implies that $\overline{v}(x) \rightarrow 0$,
regardless of the values of $a,b,c,d$. 
In case $L = \infty$, 
$v(x) \rightarrow \infty$ implies that $\overline{v}(x) \rightarrow
\infty$,
for any values of $a,b,c,d$ satisfying $a^2 +c^2 = b^2 +d^2  >
|ab+cd|$.
Thus the condition (\ref{L}), with $L =0$ or $L = \infty$, does not
characterize the pair $y_1,y_2$ up to orthogonal transformation.
In as much as this excludes some important examples, such as the
Airy equation ($q(x) =x$) and the Cauchy--Euler equation ($q(x) =
\gamma^2/x^2$), we now present a theorem which characterizes the
pair $y_1,y_2$ by a condition on $v'$:

{\bf Theorem 1.} {\em Let (\ref{q}) be an oscillatory equation and
let $y_1,y_2$ be a pair of linearly independent solutions
normalized by the unit Wronskian $|w(y_1,y_2)| =1$.  Let $v(x) =
y_1^2(x) + y_2^2(x)$ and suppose that 
\begin{equation}
\lim_{x \rightarrow \infty} v'(x) = K,\;\;0 \le K < \infty.
\label{42}
\end{equation}
Then the pair
$y_1,y_2$ is unique up to orthogonal transformation.}

{\em Proof.} Let $\alpha$ denote a phase of equation
(\ref{q})
corresponding to $y_1$ and $y_2$.  In accordance with \S 2, we have
$\alpha
\in C^3[x_0,\infty),\;\alpha'(x) \ne 0$ on $[x_0,\infty)$, and
\begin{equation}
\lim_{x \rightarrow \infty} |\alpha(x)| = \infty, \label{6}
\end{equation}
because (\ref{q}) is an oscillatory equation.  Moreover,
$$ y_1(x) = \varepsilon | \alpha'(x)|^{-1/2}\sin \alpha(x), $$
$$ y_2(x) = \varepsilon | \alpha'(x)|^{-1/2}\cos \alpha(x), $$
$\varepsilon = \pm 1$.  Hence $v(x) = y_1^2(x) + y_2^2(x) =
|\alpha'(x) |^{-1} $.
Now, for a pair $\overline{y}_1, \overline{y}_2 $, we get, from
(\ref{star})
\begin{eqnarray*}
\overline{v}'(x) &=& 
v'(x) \left[ (a^2 +c^2) \sin^2 \alpha (x) +2(ab+cd) \sin \alpha(x)
\cos
\alpha(x) + (b^2 +d^2) \cos^2 \alpha(x) \right] 
\\ & & + \frac{1}{|\alpha'(x)|}
\left[ (a^2 +c^2 -b^2 -d^2) \sin 2 \alpha (x) + 2(ab+cd) \cos
2\alpha(x) \right] \cdot \alpha'(x)
\end{eqnarray*}
i.e., 
\begin{equation}
\overline{v}'(x) = v'(x) F(x) + G(x) \label{imp}
\end{equation}
where 
$$ F(x) = (a^2 +c^2) \sin^2 \alpha (x) +2(ab+cd) \sin \alpha(x)
\cos \alpha(x) + (b^2 +d^2) \cos^2 \alpha(x) $$
is a bounded function on $[x_0, \infty)$, and 
$$ G(x) = K_1 \sin 2 \alpha (x) + K_2 \cos 2\alpha(x), $$ where
$K_1 = a^2 +c^2 -b^2 -d^2$, $K_2 = 2(ab +cd)$.
Now we show that if 
\begin{equation}
\lim_{x \rightarrow \infty} {\overline v}'(x) = {\overline
K},\;\;|{\overline K}| < \infty. \label{4222}
\end{equation}
i.e., if 
(\ref{42}) holds for $\overline{y}_1,
\overline{y}_2 $, then   $K_1 = K_2 =0$.

We have
$$
2F(x) = (a^2 +b^2 +c^2 +d^2) - K_1 \cos 2\alpha (x)
+ K_2 \sin 2\alpha (x).$$  Then, using (\ref{imp}), we get
\begin{equation}2\overline{v}'(x) = v'(x)[a^2 +b^2 +c^2 +d^2] + 
\cos 2\alpha (x)[K_2 - v'(x)K_1] +
\sin 2\alpha (x)[K_1 + v'(x)K_2].
\end{equation}
Letting $x \rightarrow \infty$, we get $K_1 = K_2 =0$. 
  Hence $a^2 +c^2 = b^2 +d^2 = K_3 \ne 0$, and $ab+cd=0$, so
$$ \left( \begin{array}{cc} a &b \\c&d \end{array} \right) \cdot
\left( \begin{array}{cc} a &c \\b&d \end{array} \right) =
\left( \begin{array}{cc} K_3 &0 \\ 0 &K_3 \end{array} \right).$$
From the normalization of the Wronskian, we get
$|ad-bc| =K_3 =1$; hence $(y_1,y_2)$ is unique up to orthogonal
transformation. This completes the proof of Theorem 1.

\subsection{An extended definition}
{\bf Definition.} {\em The  pair $y_1,y_2$ will be called principal
for an oscillatory equation (\ref{q}) 
if $v(x) := y_1^2(x) +y_2^2(x)$
satisfies either 
\begin{equation}
\lim_{x \rightarrow \infty} v(x) = L,\;(0 < L < \infty),
\label{con1}
\end{equation}
or
\begin{equation}
\lim_{x \rightarrow \infty} v'(x) = K,\;(0 \le K < \infty).
\label{con4}
\end{equation}
}
  Due to \cite{ENV}, and Theorem 1, this
pair is unique up to  orthogonal transformation.

Since all principal pairs for an oscillatory equation (if they
exist) differ only by an orthogonal transformation, the expression
$$v(x) = y_1^2(x) +y_2^2(x)$$
remains the same for all principal pairs.

The following sufficient condition for the existence of principal
pairs is simply a restatement of a result of P. Hartman
\cite[Theorem $20.1_0$]{H61}:

{\bf Corollary 1.} {\em Let $q \in C^1[x_0,\infty)$, $q'(x) \ge
0$, $q''(x) \le 0$ and $\lim_{x \rightarrow \infty} q(x) =
\infty$. Then (\ref{q}) has a principal pair of solutions $y_1,y_2$
such that 
$$v(x) \ge 0,\;\;\lim_{x \rightarrow \infty}v(x) =0,$$
$$v'(x) \le 0,\;\;\lim_{x \rightarrow \infty}v'(x) =0,$$
$$\lim_{x \rightarrow \infty}v''(x) =0 $$ and
$$\lim_{x \rightarrow \infty}v'''(x) =0. $$}

Another sufficient condition for the existence of principal
pairs follows from another  result of Hartman
\cite[Theorem $22.1_0$]{H61}:

{\bf Corollary 2.} {\em Let $q \in C^2[x_0,\infty)$, $q'(x) \le
0$, $q(x)q''(x) -3q'(x)^2 \ge 0$.  Then (\ref{q}) has a principal
pair of solutions $y_1,y_2$ such that 
$$v(x) \ge 0,$$
$$\lim_{x \rightarrow \infty} \left[ q^{-1}(x) \left(\frac{d}{dx}
\right)^nv(x)\right] = 0,\;n=1,2,3,$$ and $v$ approaches a finite
limit or $\infty$ according as $q$ approaches a positive limit or
$0$. }

\subsection{Examples}
{\bf Example 1.} The generalized Airy equation

This is an example with $L = 0$. The equation in question is 
\begin{equation}
y'' + (2\nu)^{-2}x^{1/\nu -2}y =0, \label{gairy}
\end{equation}
where we suppose that $0 < \nu \le 1/2$. The usual Airy equation
corresponds to $\nu = 1/3$. With the usual notation for the Bessel
functions \cite{wat}, the pair  $x^{1/2}J_\nu(2\nu x^{1/(2\nu)}),\;
x^{1/2}Y_\nu(2\nu x^{1/(2\nu)})$ satisfy (\ref{gairy}).
In view of the known result \cite[p. 446]{wat} that, for $0 \le \mu
< 1/2$,  the function $t[J_\mu^2(t)
+ Y_\mu^2(t)]$ increases to $2/\pi$ on $(0,\infty)$, we see that
$$v(x) = x[J_\nu^2(2\nu x^{1/(2\nu)}) + Y_\nu^2(2\nu
x^{1/(2\nu)})]$$ approaches $0$ as $x \rightarrow \infty$.  In
fact, for
$\frac13 \le \nu < \frac12$,  $v(x)$ is
completely monotonic (i.e., its successive derivatives alternate in
sign) \cite[Theorem 5.1 (i)]{muldoon}. The asymptotic formula
\begin{equation}
x[J_\nu^2(2\nu x^{1/(2\nu)}) + Y_\nu^2(2\nu x^{1/(2\nu)})] =
\frac{1}{\nu\pi}x^{1-1/(2\nu)}\left[ 1 + O(x^{-1/\nu}) \right],
\end{equation} \cite[p. 224 (5)]{wat} shows that $v(x)$ and $v'(x)$
approach
$0$ as $x \rightarrow \infty$. Thus, by Theorem 1, 
the pair $x^{1/2}J_\nu(2\nu x^{1/(2\nu)}), x^{1/2}Y_\nu(2\nu
x^{1/(2\nu)})$ constitute a principal pair for (\ref{gairy}).

{\bf Example 2.} 

The equation
\begin{equation}
y'' + x^{-1} y =0, \label{j1}
\end{equation}
has a principal pair  $x^{1/2}J_1(2 x^{1/2}), x^{1/2}Y_1(2
x^{1/2})$.  Here $v(x) \rightarrow \infty$, but $v'(x) \rightarrow
0$ as $x \rightarrow \infty$, so the hypotheses of Theorem 1 are
satisfied again.

{\bf Example 3.} The Cauchy--Euler equation

This is an example with $L = \infty$, to which Theorem 1 is 
applicable. The Cauchy--Euler equation is
\begin{equation}
y'' + \frac{\gamma^2}{x^2}y =0. \label{cee}
\end{equation}
In case $\gamma^2 > 1/4$, the equation is oscillatory and a 
principal pair of solutions is given by $y_1(x)= x^{1/2}\sin(s\log
x)$, $y_2(x)=
x^{1/2}\cos(s\log x)$, where $s = \sqrt{\gamma^2 -1/4}$, so $v(x)
= x$ and $v'(x) = 1$, so the hypotheses of Theorem 1 are satisfied.

It is clear that every pair 
$$ \overline{y}_1 = ay_1 + by_2,$$
$$ \overline{y}_2 = cy_1 + dy_2,$$
with 
$a^2 +c^2 = b^2 +d^2  > |ab+cd|$ has the property that
$\overline{v}(x) \rightarrow \infty$ as $x \rightarrow \infty$, but
it is only in the case $a^2 +c^2 = c^2 +d^2,\;ab+cd =0$ that 
$\overline{v}'(x)$ approaches a finite limit.

\section{Properties of principal pairs}
Suppose that an oscillatory equation (1)
admits a pair
$y_1, y_2$, for which
\begin{equation}
\lim _{x \rightarrow \infty} v'(x) = \lim _{x \rightarrow \infty}
(y_1^2(x) + y_2^2(x))' = 0.
\label{81}
\end{equation}
Due to the extended definition in \S 3.1, the pair $ y_1, y_2 $
is a principal pair.
Let $\alpha$ denote the phase
of (1) corresponding to $y_1, y_2$. Then
$$ y_1(x) = \varepsilon | \alpha'(x)|^{-1/2}\sin \alpha(x), $$
 $$ y_2(x) = \varepsilon | \alpha'(x)|^{-1/2}\cos \alpha(x), $$
$$\varepsilon = \pm 1,$$
$$\alpha
\in C^3[x_0,\infty),\;\alpha'(x) \ne 0\;{\rm on}\;[x_0,\infty),$$
 $$v(x) = y_1^2(x) +y_2^2(x)= |\alpha'(x)|^{-1}.$$
Condition (\ref{81}) gives
\begin{equation}
\lim _{x \rightarrow \infty} \frac{\alpha''(x)}{(\alpha'(x))^2} =0.
\label{222}
\end{equation}
Consider the first derivative of $y_1$:
$$
y_1'(x) = \varepsilon \left[ -\frac12 |\alpha'(x)|^{-1/2} \frac
{\alpha''(x)}{\alpha'(x)} \sin \alpha(x) +
|\alpha'(x)|^{1/2} \cos \alpha(x) \cdot \;{\rm sign}\; \alpha'
\right].
$$
All of its
zeros are at points $x_j$ where \begin{equation}
\frac{\cos \alpha(x_j)}{\sin \alpha(x_j)} =\frac12
\frac{\alpha''(x_j)}{(\alpha'(x_j))^2}. \label{11}
\end{equation}
From (\ref{222}) we have
$$|\alpha(x_j)  + \pi/2 - j\pi\; {\rm sign}\; \alpha'| \rightarrow
0, \;\;{\rm as} \; \; j \rightarrow \infty.$$
A general solution of (\ref{q}) can be written in the form
$$ y(x;c_1,c_2) = c_1| \alpha'(x)|^{-1/2}\sin (\alpha(x) +c_2). $$
Its zeros are at points $x_j$, where 
$$
\alpha(x_j) +c_2 = j\pi \;{\rm sign}\;\alpha',\;j=j_0,j_0 +1,\dots.
$$
Hence
$$ y(x;c_1,\pi/2) = c_1| \alpha'(x)|^{-1/2}\cos \alpha(x) $$
is that solution (up to constant multiple) whose zeros are
approached by  the
zeros of the derivative of $y_1$. However, after the normalization
$c_1 = \pm 1$, $y_1$ together with $y(x, \pm 1, \pi/2) = \pm
y_2(x)$
represent a principal pair of (\ref{q}). 

We may summarize our considerations in the following Theorem:

{\bf Theorem 2.} {\em If an oscillatory equation (\ref{q}) has a
pair of solutions $y_1,y_2$ for which
\begin{equation}
\lim_{x \rightarrow \infty} v'(x) = 0, \label{42424}
\end{equation}
then they have the property that the difference of the zeros of the
derivative of one
of them and the zeros of the other tend to $0$ as $ x \rightarrow
\infty$.  In other words, among all sequences of zeros of solutions
of
(\ref{q}), the
best approximation to the zeros of $y_1'$ is given by zeros of
$y_2$. }

{\em Remark.} Evidently, the conditions in Corollary 1 are
sufficient to guarantee the conclusion of Theorem 2. Further
sufficient conditions  can be obtained by using Hartman
\cite[Theorem $18.1_0$]{H61}.  In fact, this Theorem,
together with his Theorem $20.1_0$, show that the
requirements
$$ q' \ge 0,\;q'' \le 0,\;{\rm and}\; \lim _{x \rightarrow \infty}
q(x) = q(\infty),\;
{\rm where}\; 0 < q(\infty) < \infty$$
 are sufficient to imply the conclusion of our Theorem 2.

The property described in Theorem 2 links two transformations of
second order equations.

 The first is the Kummer (sometimes called St\"ackel)
transformation consisting in the change of independent and
dependent variables of the solution $z(t)$ into the solution $y(x)$
expressible in the form
$$y(x) = f(x) z(h(x))$$
or, better,
$$y_1(x) = f(x) z_1(h(x))$$ and
$$y_2(x) = f(x) z_2(h(x)),$$
for two linearly independent solutions $y_1,y_2$ and $z_1,z_2$ of
the transformed equations.  In particular,
$$y_1(x) = f(x) \sin h(x)$$ and
$$y_2(x) = f(x) \cos h(x),$$
when the equation $z''+z=0$ is taken as canonical (as can always be
done). 
The second is the Pr\"ufer transformation, considering again a two
dimensional underlying space for two functions, now 
$$y_1(x) = y(x) = \rho(x) \sin \varphi(x)$$
and
$$y_2(x) = y'(x) = \rho(x) \cos \varphi(x).$$

\end{document}